\documentclass[12pt]{amsart} 
\usepackage{amssymb,amsmath,latexsym, enumerate} 

\usepackage{color}

\numberwithin{equation}{section}
\theoremstyle{plain} 
\newtheorem{thm}{Theorem}[section]

\newtheorem{lem}[thm]{Lemma}
\newtheorem{rmk}[thm]{Remark}
\newtheorem{ex}[thm]{Example}
\newtheorem{cor}[thm]{Corollary}
\newtheorem{prop}[thm]{Proposition}

\newtheorem{df}[thm]{Definition}

\begin{document}
\newcommand{\hgt}{\operatorname{ht}}
\newcommand{\ext}{\operatorname{Ext}}
\newcommand{\Hom}{\operatorname{Hom}}
\newcommand{\extx}{\operatorname{Ext}_{R_x}}

\newcommand{\exti}{\operatorname{Ext}^i_R}
\newcommand{\extijM}{\operatorname{Ext}^i_R(R/I,H^j_I(M))}
\newcommand{\extjiM}{\operatorname{Ext}^j_R(R/I,H^i_I(M))}

\newcommand{\extijiM}{\operatorname{Ext}^i_R(R/I,H^{j+1}_I(M))}
\newcommand{\extijR}{\operatorname{Ext}^i_R(R/I,H^j_I(R))}
\newcommand{\extidR}{\operatorname{Ext}^i_R(R/I,H^d_I(R))}
\newcommand{\extidiR}{\operatorname{Ext}^i_R(R/I,H^{d-1}_I(R))}

\newcommand{\hjR}{H^j_I(R)}
\newcommand{\hiR}{H^i_I(R)}
\newcommand{\hjM}{H^j_I(M)}
\newcommand{\hjiM}{H^{j+1}_I(M)}
\newcommand{\hiM}{H^i_I(M)}
\newcommand{\hhM}{H^h_I(M)}
\newcommand{\hnM}{H^n_I(M)}
\newcommand{\hniM}{H^{n-1}_I(M)}

\newcommand{\hjRx}{H^j_I{R_x}(R_x)}

\newcommand{\hdM}{H^d_I(M)}
\newcommand{\hdR}{H^d_I(R)}
\newcommand{\hdiM}{H^{d-1}_I(M)}
\newcommand{\hdiR}{H^{d-1}_I(R)}
\newcommand{\hdiiR}{H^{d-2}_I(R)}
\newcommand{\hiIM}{H^i_I(M)}
\newcommand{\hgIM}{H^g_I(M)}
\newcommand{\supp}{\operatorname{Supp}_R}
\newcommand{\ass}{\operatorname{Ass}_R}
\newcommand{\suppx}{\operatorname{Supp}_{R_x}}
\newcommand{\assx}{\operatorname{Ass}_{R_x}}
\newcommand{\depth}{\operatorname{depth}}
\newcommand{\ch}{\operatorname{v}}
\newcommand{\ann}{\operatorname{ann}}

\title[Cofiniteness and associated primes]{Cofiniteness and associated primes of local cohomology modules}

\author{Thomas Marley}
\address{University of Nebraska-Lincoln\\
Department of Mathematics and Statistics\\
Lincoln, NE 68588-0323}
\email{tmarley@math.unl.edu}
\urladdr{http://www.math.unl.edu/\textasciitilde tmarley}
\author{Janet C. Vassilev}
\address{University of Arkansas\\
Department of Mathematical Sciences\\
Fayetteville, AR 72701}
\email{jvassil@uark.edu}
\urladdr{http://comp.uark.edu/\textasciitilde jvassil}

\thanks{The first author was partially supported by NSF grant DMS-0071008.}
%\keywords{Local Cohomology, cofinite, associated prime}
%\subjclass[2000]{13D45}
\date{January 15, 2002}

\begin{abstract}  Let $R$ be a $d$-dimensional regular local ring, $I$ an ideal of
$R$, and $M$ a finitely generated $R$-module of dimension $n$.  We prove that the set of
associated primes of $\ext^i_R(R/I,H^j_I(M))$ is finite for all $i$ and $j$ in the following
cases:
\begin{enumerate}[(a)]
\item $\dim M\le 3$;
\item $\dim R\le 4$;
\item $\dim M/IM\le 2$ and $M$ satisfies Serre's condition $S_{n-3}$;
\item $\dim M/IM\le 3$, $\operatorname{ann}_RM=0$, $R$ is unramified, and $M$ satisfies $S_{n-3}$.
\end{enumerate}
In these cases we also prove that $H^i_I(M)_p$ is $I_p$-cofinite for all
but finitely many primes $p$ of $R$. 

Additionally, we show that if $\dim R/I\ge 2$ and
$\operatorname{Spec}R/I-\{m/I\}$ is disconnected then
$H^{d-1}_I(R)$ is not $I$-cofinite, generalizing a result due to Huneke and Koh.   
\end{abstract}

\maketitle

\section{Introduction}

Let $R$ be a local Noetherian ring, $I$ an ideal of $R$ and $M$ a finitely generated 
$R$-module.  It is well known that the local cohomology
modules $\hiM$ are not generally Noetherian for $i>0$.  However, in recent
years there have been several results showing that, under certain conditions, local cohomology
modules share some special properties with Noetherian modules.   
Perhaps the most striking of these results is the following:  If $R$ is an unramified regular
local ring then for all ideals $I$ of $R$ and
all $i\ge 0$ the set of associated primes of $H^i_I(R)$ is finite and the Bass numbers (with respect
to any prime) of $H^i_I(R)$ are finite.  This was proved in the case
of positive characteristic by Huneke and Sharp \cite{HS} and 
in the case of characteristic zero
by Lyubeznik \cite{L1}, \cite{L2}.  

In a 1970 paper Hartshorne \cite{Ha} gave an example which 
shows that the Bass numbers of $H^i_I(R)$ may be infinite if $R$ is not regular.  Until
recently, it was
an open question as to whether the set of associated primes of $H^i_I(R)$ is finite
for an arbitrary Noetherian ring $R$ and any ideal $I$.  However, examples given  
by A. Singh \cite{Si} (in the non-local case) and M. Katzman \cite{Kz} (in the local case)
show there exist local cohomology modules of
Noetherian rings with infinitely many associated primes.

In this paper, we show that when $\dim R$ or $\dim M/IM$ is `small' the set of
associated primes of $H^i_I(M)$ for $i\ge 0$ is
finite.  In fact, 
we are able to prove the stronger result that the set
of associated primes of $\ext^i_R(R/I,H^j_I(M))$ is finite for all $i$ and $j$.  Furthermore,
we show that $H^j_I(M)$ is `almost' locally $I$-cofinite, meaning that for all but finitely
many primes $p$ of $R$, $H^j_I(M)_p$ is $I_p$-cofinite.
We can summarize many of our main results (found in section 3) as follows:

\begin{thm} Let $R$ be a regular local ring, 
$I$ an ideal of $R$, and $M$ a finitely
generated $R$-module of dimension $n$.  Suppose one of the following conditions holds:
\begin{enumerate}[(a)]
\item $\dim M\le 3$;
\item $\dim R\le 4$;
\item $\dim M/IM\le 2$ and $M$ satisfies Serre's condition $S_{n-3}$;
\item $\dim R/I=3$, $\operatorname{ann}_RM=0$, $R$ is unramified, and $M$ satisfies $S_{n-3}$.
\end{enumerate}
Then $\ext^i_R(R/I,H^j_I(M))$ has finitely many associated primes for all $i$ and $j$.
Furthermore, $H^j_I(M)_p$ is $I_p$-cofinite for almost all primes $p$; i.e.,
$\ext^i_R(R/I,H^j_I(M))_p$ is a finitely generated $R_p$-module for all
$i,j$ and all but finitely many $p$.
\end{thm}

We give examples in section 3 to show that some of the above conditions for $H^j_I(M)$ to
be almost locally $I$-cofinite are the best possible.
Also, Katzman's example shows that the finiteness of the set of associated primes of $H^j_I(R)$
does not hold for arbitrary modules over a six-dimensional regular local ring.
It remains an open question whether there exists a local cohomology module over
a five-dimensional regular local ring with infinitely many associated primes.

Of central importance in this paper is the concept of cofiniteness, first defined by Hartshorne
\cite{Ha} and further studied by Huneke and Koh \cite{HK}. 
Many of our ideas were inspired
by a careful reading of these two papers.  An $R$-module $N$ is called $I$-cofinite if
$\supp N\subseteq \operatorname{V}(I)$ and $\ext^i_R(R/I,N)$ is finitely generated for
all $i\ge 0$. (Here and throughout, $\operatorname{V}(I)$ denotes
the set of prime ideals containing $I$.) Hartshorne proved that if $R$
is complete regular local ring and $M$ is a finitely generated $R$-module then
the local cohomology modules $H^i_I(M)$ are $I$-cofinite if either $\dim R/I\le 1$ or $I$ is
principal. Subsequently, these results were shown 
to hold for arbitrary Noetherian local rings
in \cite{Yo}, \cite{DM}, and \cite{Ka}.  
Hartshorne, and later Huneke and Koh, gave examples
of local cohomology modules $H^i_I(M)$ which are not $I$-cofinite.  
A key step in many of our proofs is to
choose an element $x\in R$ avoiding a countably infinite set of primes such that $H^i_I(M)_x$
is $I_x$-cofinite for all $i$.  When such an element exists is the subject
of section 2.

In section 3 we prove Theorem 1.1 as well as  the following generalization of
\cite[Theorem 3.6(ii)]{HK}:

\begin{thm} Let $(R,m)$ be a complete Cohen-Macaulay normal local ring and $I$ an ideal
such that $\dim R/I\ge 2$ and $\operatorname{Spec}R/I-\{m/I\}$ is disconnected.  Then
$\operatorname{Hom}_R(R/I, H^{d-1}_I(R))$
is not finitely generated.  Consequently, $H^{d-1}_I(R)$ is not $I$-cofinite.
\end{thm}

By a `ring' we always mean a commutative ring with identity.   Local rings are assumed to be
Noetherian.  We refer the reader to \cite{Mat} or \cite{BH}
for any unexplained terminology.  In particular, $\ass M$ denotes the set of associated primes
of the $R$-module $M$ and $\supp M$ denotes the support of $M$. 
For the definition of local cohomology and its basic properties, we refer the reader
to \cite{BS}.

\section{Preliminary Results}

We begin this section with an elementary result:

\begin{lem} \label{ideal} Let $R$ be a Noetherian ring and $M$ an $R$-module.  Then the set
$$\{x\in R\mid M_x \text{ is a finitely generated }R_x\text{-module}\}$$
is an ideal of $R$.
\end{lem}

{\it Proof:} It is enough to show that if $M_x$ and $M_y$ are finitely generated modules 
(over $R_x$ and $R_y$, respectively), then $M_{x+y}$ is a finitely generated $R_{x+y}$-module.
Let $A,B$ be finitely generated $R$-submodules of $M$ such that $A_x=M_x$ and $B_y=M_y$.
We claim that $(A+B)_{x+y}=M_{x+y}$:  Let $m\in M$.  Since $R_xm\subseteq A_x$ we have $x\in
\sqrt{(A:m)}$, where $(A:m)=\{r\in R\mid rm\in A\}$.  Similarly, $y\in \sqrt{(B:m)}$.
Hence, $x+y\in \sqrt{(A+B:m)}$, which implies $R_{x+y}m\subseteq (A+B)_{x+y}$.
\qed

\begin{df}{\rm Let $R$ be local ring of
dimension $d$, $I$ an ideal, and $M$ a finitely generated $R$-module.
For each $i\ge 0$ define
$$D^i(I,M):=\{x\in R\mid H^i_I(M)_x \text{ is }I_x\text{-cofinite}\}.$$
By the lemma, it is clear that $D^i(I,M)$ is an ideal of $R$.  Further, let
$$D(I,M):=\bigcap_{i=0}^dD^i(I,M).$$}
\end{df}

We note that if $p\not\supset D^i=D^i(I,M)$ then $H^i_I(M)_p$ is $I_p$-cofinite.
Thus, $\operatorname{V}(D^i)$ contains the non-$I$-cofinite locus of $H^i_I(M)$.
The following are some additional observations concerning $D^i(I,M)$.
Throughout this section, we adopt the convention
that the dimension of the zero ring is $-1$. 

\begin{rmk} \label{D-rmk}  Let $R$ be a Noetherian local ring, $I$ an ideal of $R$,
and $M$ a finitely generated $R$-module of dimension $n$.  Let $D^i=D^i(I,M)$.  Then
\begin{enumerate}[(a)]
\item $H^i_I(M)$ is $I$-cofinite if and only if $D^i=R$;
\item $D^i$ is a radical ideal containing $I$;
\item $D^0=R$;
\item $\dim R/D^n\le 0$;
\item $\dim R/D^{n-1}\le 1$.
\end{enumerate}
\end{rmk}

{\it Proof:} Statements (a)-(d) are clear.  For (e), we note that $\supp H^{n-1}_I(M)$ is finite by \cite[Corollary 2.5]{Mar}.
\qed

The main result of this section is that under certain conditions
$\dim R/D\le 1$.  In section 3 we give examples to show that this inequality does not
hold in general. 

\begin{thm}\label{bigprop} Let $R$ be a local ring, $I$ an ideal of $R$, and
$M$ a finitely generated $R$-module of dimension $n$.  Let $D:=D(I,M)$ and suppose one of the following conditions holds:
\begin{enumerate}[(a)]
\item $\dim M\le 3$; 
\item $\dim R=4$ and $R$ is a UFD; 
\item $R$ is the quotient of a Cohen-Macaulay ring, $\dim M/IM\le 2$, and either $\dim M\le 4$ or 
$M$ satisfies Serre's condition $S_{n-3}$; 
\item $R$ is a unramified regular local ring, $\dim R/I\le 3$, and $M$ satisfies $S_{d-3}$ where
$d=\dim R=\dim M$.
\end{enumerate}
Then $\dim R/D\le 1$.
\end{thm}

The proof of case (a) follows immediately from parts (c), (d), and (e) of Remark \ref{D-rmk}. 
We give the proofs of the remaining cases separately. 
The following proposition will be useful in our arguments.

\begin{prop} \label{Icof}  Let $R$ be a Noetherian ring, $I$ an ideal of $R$, and $M$ a
finitely generated $R$-module.  Suppose there exists an integer $h\ge 0$ such that
$\hiM$ is $I$-cofinite for all $i\neq h$.
Then $\hiM$ is $I$-cofinite for all $i$.
\end{prop}

{\it Proof:} 
Consider the Grothendieck spectral sequence
$$E_2^{p,q}=
\operatorname{Ext}^p_R(R/I, H^q_I(M)) \Rightarrow \operatorname{Ext}^{p+q}_R(R/I,M).$$ 
Since $E_r^{p,q}$ is a subquotient of $E_2^{p,q}$ for all $r\ge 2$,
our hypotheses give us that $E_r^{pq}$ is finitely generated for all $r\ge 2$, $p\ge 0$,
and $q\neq h$.  For each $r\ge 2$ and $p,q\ge 0$, let $Z_r^{p,q}
=\ker (E_r^{p,q}\to E_r^{p+r,q-r+1})$
and $B_r^{p,q}=\operatorname{im} (E_r^{p-r,q+r-1}\to E_r^{p,q})$.   Note that $B_r^{p,q}$
is finitely generated for all $p$, $q$, and $r\ge 2$, since either $E_r^{p-r,q+r-1}$ or
$E_r^{p,q}$ is finitely generated.
For all $r\ge 2$ and 
$p\ge 0$
we have the exact sequences
$$0\rightarrow B_r^{p,h}\rightarrow Z_r^{p,h}\rightarrow E_{r+1}^{p,h}\rightarrow 0$$
and
$$0\rightarrow Z_r^{p,h}\rightarrow E_r^{p,h}\rightarrow B_r^{p+r,q-r+1}\rightarrow 0.$$
Now $E_{\infty}^{p,h}$ is isomorphic to a subquotient of $\operatorname{Ext}^{p+h}_R(R/I,M)$
and thus is finitely generated for all $p$.  
Since $E_r^{p,h}=E_{\infty}^{p,h}$ for $r$ sufficiently large, we have that $E_r^{p,h}$ is
finitely generated for all $p$ and all large $r$. Fix $p$ and $r$ and
suppose $E_{r+1}^{p,h}$ is finitely generated.
>From the first exact sequence we obtain that $Z_r^{p,h}$ is finitely generated. 
 From the second exact sequence we get that $E_r^{p,h}$ is finitely generated.  Continuing in
this fashion, we see that $E_r^{p,h}$ is finitely generated for all $r\ge 2$ and all $p$.  
In particular, $E_2^{p,h}=\operatorname{Ext}^p_R(R/I,H^h_I(M))$ is finitely generated for
all $p$. 
\qed

We note some easy consequences of this result:

\begin{cor} \label{Icofcor} Let $R$ be a Noetherian ring, $I$ an ideal of $R$, $M$ a
finitely generated $R$-module, and $h\in \mathbb Z$. 
\begin{enumerate}[(a)]
\item If $R$ is local and $\dim R\le 2$ then $H^i_I(M)$ is $I$-cofinite for all $i\ge 0$.
\item Suppose 
$\hiM$ is finitely generated for all $i<h$ and $\hiM=0$
for all $i>h$.  Then $\hiM$ is $I$-cofinite for all $i$.  In particular, if $H^i_I(M)=0$ for
all $i>1$ then $\hiM$ is $I$-cofinite for all $i$ (cf. \cite{Ha}, \cite{Ka}).
\item Suppose $\hiM=0$ for all $i\neq h,h+1$.  Then $H^h_I(M)$ is $I$-cofinite
if and only if $H^{h+1}_I(M)$ is $I$-cofinite.
\end{enumerate}
\end{cor}

{\it Proof:} For part (a), note that $H^0_I(M)$ and $H^d_I(M)$, where $d=\dim R$, are both
$I$-cofinite (\cite[Theorem 3]{DM}).  The remaining statements follow immediately from
Proposition \ref{Icof}.
\qed

We will also need the following result:

\begin{lem} \label{fg} Let $R$ be a local ring 
which is the homomorphic image of a Cohen-Macaulay ring,
$I$ an ideal of $R$, and $M$ a finitely generated $R$-module.  Let $n=\dim M$ and $r=\dim M/IM$.
Suppose that 
\begin{enumerate}[a)]
\item $M$ is equidimensional, and
\item $M$ satisfies Serre's condition $S_{\ell}$ for some $\ell\le n-r-1$.
\end{enumerate}
Then $H^i_I(M)$ is finitely generated for all $i<\ell+1$.
\end{lem}

{\it Proof:} Without loss of generality we can assume that $\operatorname{ann}_RM=(0)$.  As $M$ is an equidimensional $R$-module with $\operatorname{ann}_RM=(0)$, we note
$R$ is equidimensional.  As $R$ is a quotient
of a Cohen-Macaulay ring, we may assume $R$ is complete as well (\cite[9.6.3]{BS}, 
\cite[\S 23, \S 31]{Mat}).  Let $p\in \operatorname{Spec}R$.
If $\hgt p\le \ell$ then 
\begin{align*}
\depth M_p + \hgt (I+p)/p&=\hgt p +\hgt (I+p)/p\\
&=\hgt (I+p)\\
&\ge \hgt I \\
&=n-r\\
&\ge \ell+1.
\end{align*}
If $\hgt p>\ell$ and $p\not\supset I$ then 
$\depth M_p+\hgt (I+p)/p\ge \ell + 1.$
Hence
$$\min \{\depth M_p + \hgt (I+p)/p\mid p\not\in \operatorname{V}(I)\}\ge \ell+1.$$
By \cite{Fa} or \cite[9.5.2]{BS}, $H^i_I(M)$ is finitely generated for all $i<\ell+1$.
\qed

The following is a generalization of \cite[5.10.9]{EGA}:

\begin{lem} \label{equid} Let $(R,m)$ be a catenary local ring and $M$ a
finitely generated indecomposable $R$-module
which satisfies Serre's condition ${\rm S}_2$.  Then $M$ is equidimensional.
\end{lem}

{\it Proof:} Without loss of
generality, we may assume $\operatorname{ann}_R M=0$.  Suppose $R$ is not equidimensional.
Let $X:=\{P\in \operatorname{Min} R\mid \dim R/P=\dim R\}$ and
$Y:=\{P\in \operatorname{Min}R \mid \dim R/P<\dim R\}$.  Furthermore, let
$$I=\bigcap_{P\in X} P \qquad \text{ and } \qquad
J=\bigcap_{P\in Y} P.$$
Then $I+J\subseteq m$ and $\hgt (I+J)\ge 1$.  

We claim that $\hgt (I+J) \ge 2$.
For, suppose there exists a height one prime $Q$ containing $I+J$. Then there exist
$P_1\in X$ and $P_2\in Y$ such that $P_1+P_2\subseteq Q$.
 As $R$ is local, catenary and $\hgt Q/P_1=\hgt Q/P_2=1$, we have
$$\dim R/P_1=\dim R/Q + \hgt Q/P_1 =\dim R/Q + \hgt Q/P_2=\dim R/P_2,$$
a contradiction. 

Since $\hgt (I+J)\ge 2$ and $M$ satisfies $S_2$, we have $H^0_{I+J}(M)=H^1_{I+J}(M)=0$.  Thus,
by the Mayer-Vietoris sequence, we obtain
$$M=H^0_{I\cap J}(M)\cong H^0_I(M)\oplus H^0_J(M).$$
But $H^0_I(M)$ and $H^0_J(M)$ are nonzero, contradicting that $M$ is indecomposable.
\qed

The condition that $M$ be indecomposable in the above lemma is necessary.  For example, let 
$R=k[[u,v,w,x,y]]$ be a power series ring over a field $k$.  Then
$M=R/(u,v)\oplus R/(w,x,y)$ is $S_2$ but not equidimensional.

We now give the proof of part (c) of Theorem \ref{bigprop}.

\begin{prop} \label{bigprop3} Let $R$ be a local ring  
which is the homomorphic
image of a Cohen-Macaulay ring, 
$M$ a finitely generated $R$-module of dimension $n$, and $I$ an ideal of $R$ such that
$\dim M/IM\le 2$.  Suppose that either $n\le 4$ or $M$ satisfies $S_{n-3}$.  Then 
$\dim R/D\le 1$ where $D=D(I,M)$.
\end{prop}

{\it Proof:}  We consider first the case when $n\le 4$.  If $n\le 3$ the result follows by
Theorem \ref{bigprop} (a).
Assume that $\dim R/p=4$ for all $p\in \ass M$.  Then $M$ satisfies $S_1$ and $H^i_I(M)$
is finitely generated for $i<2$ by Lemma \ref{fg}.  
Let $J=D^{3}(I,M)\cap D^4(I,M)$.  Then $\dim R/J\le 1$ by Remark \ref{D-rmk}.
For any $x\in J$, $\hiM_x$  
is $I_x$-cofinite for $i\neq 2$. 
By Corollary  \ref{Icofcor}, $\hiM_x$ is
$I_x$-cofinite for all $i$ and for all $x\in J$.  Thus, $J\subseteq D$ and $\dim R/D\le 1$.

Now let $M$ be an arbitrary finitely generated 
four-dimensional $R$-module.  Let $N$ be the largest $R$-submodule of $M$ such that
$\dim N\le 3$.  Then $\dim R/p=4$ for all $p\in \ass M/N$.  Thus, $\dim R/D(I,M/N)\le 1$
by the preceding argument.  Let $V(J_1)=\supp H^2_I(N)\cup \supp H^3_I(N)$ and  
 $J=J_1\cap D(I,M/N)$.  By Remark \ref{D-rmk}, $\dim R/J\le 1$.
For all
$x\in J$, $H^i_I(N)_x=0$ for $i>1$ and $H^i_I(M/N)_x$ is $I_x$-cofinite
for all $i$.  Hence $H^i_I(M)_x\cong H^i_I(M/N)_x$ for all $i>1$; consequently,
$H^i_I(M)_x$ is $I_x$-cofinite for all $i>1$.
 By Proposition \ref{Icof}, we obtain that $\hiM_x$ is $I_x$-cofinite
for all $i$ and all $x\in J$.  Thus, $J\subseteq D$.

We now consider the case when $n\ge 5$ and $M$ satisfies $S_{n-3}$.  
Without loss of generality, we may assume that $M$ is indecomposable.
(If $M=M_1\oplus M_2\oplus \cdots \oplus M_n$  then
$D(I,M_1)\cap \cdots \cap D(I,M_n)
=D(I,M)$.)  Since $R$ is catenary, $M$ is equidimensional by Lemma \ref{equid}.
Hence, by Lemma \ref{fg}, $\hiM$ is finitely generated for $i<n-2$.
Let $J=D^{n-1}(I,M)\cap D^n(I,M)$.  Then $\dim R/J\le 1$ and for all $x\in J$
and $i\neq n-2$ we have that $\hiM_x$ is $I_x$-cofinite.  Thus, by Proposition \ref{Icof}, 
$\hiM_x$ is $I_x$-cofinite for all $i\ge 0$ and $x\in J$.  
\qed

We next give the proof of part (b) of Theorem \ref{bigprop}:

\begin{prop} \label{bigprop2} Let $R$ be a four-dimensional local UFD, $I$ an ideal of $R$, and $M$ a
finitely generated $R$-module.  Then $\dim R/D\le 1$ where $D=D(I,M)$.
\end{prop}

{\it Proof:}  We may assume that $\operatorname{dim} \ M=4$ by Theorem \ref{bigprop} (a).  The case when $\hgt I=0$ is trivial ($D=R$).

{\it Case 1:}  $\hgt I\ge 2$.  

Suppose first that $M$ is torsion-free.  Since $R$ is a domain there exists a torsion $R$-module
$C$ and a short exact sequence
$$0 \rightarrow M \rightarrow R^n \rightarrow C \rightarrow 0,$$ where
$n=\operatorname{rank}M$.  Since $R$ satisfies Serre's condition $S_2$
we have $H^i_I(R)=0$ for $i=0,1$.  Therefore
$H^1_I(M) \cong H^0_I(C)$, which is a finitely generated $R$-module.
Let $J=D^3(I,M)\cap D^4(I,M)$.  Then $\dim R/J\le 1$ and 
$H^i_I(M)_x$ is $I_x$-cofinite for all $x\in J$ and $i\neq 2$.
Therefore, by Proposition \ref{Icof},  $\hiM_x$ is $I_x$-cofinite for all $i$ and $x\in J$.
Consequently, $J\subseteq D$.

Now let $M$ be an arbitrary finitely generated $R$-module and $T$ the
torsion submodule of $M$.  Then $M/T$ is torsion-free and $\dim T\le 3$.  Let $J=J_1\cap J_2$ where
$J_1=D(I,M/T)$ and 
$\operatorname{V}(J_2)=\supp H^2_I(T) \cup \supp H^3_I(T)$.
Then $\dim R/J\le 1$.  Let $x\in J$.  From the exact sequence 
$$0\to T_x\to M_x\to (M/T)_x\to 0$$
we obtain that $H^i_I(M)_x\cong H^i_I(M/T)_x$ for all $i\ge 2$.  Therefore, $\hiM_x$ is $I_x$-cofinite
for all $i\neq 1$. By Proposition \ref{Icof} $\hiM_x$ is $I_x$-cofinite for all $i$.   
Thus, $\dim R/D\le 1$.

{\it Case 2:} $\hgt I =1$. 
 
We may assume $I=\sqrt{I}$. As $R$ is a UFD,
$I=K \cap (f)$ where $\hgt K \geq 2$ and $\hgt (K+(f)) \geq 3$.   
Let $J=D(K,M)\cap (K+(f))$.  Then $\dim R/J\le 1$ by Case 1.  For $x\in J$ we have that
$K_x+(f)_x=R_x$ and $(R/I)_x\cong (R/K)_x\oplus (R/(f))_x$.  By the Mayer-Vietoris sequence,
we obtain that $$H^i_I(M)_x\cong H^i_K(M)_x\oplus H^i_{(f)}(M)_x$$ for all $i$.  Furthermore,
$$\ext^i_R(R/K,H^j_{(f)}(M))_x=0=\ext^i_R(R/(f),H^j_K(M))_x$$
for all $i$ and $j$, since any prime in the support of these modules contains $K_x+(f)_x=R_x$.
Combining these facts, we have for all $x\in J$ and all $i,j$
$$\ext^i_R(R/I,H^j_I(M))_x\cong \ext^i_R(R/K,H^j_K(M))_x\oplus \ext^i_R(R/(f),H^j_{(f)}(M))_x.$$
The first summand is finitely generated as $x\in D(K,M)$ and the second
is finitely generated by Corollary \ref{Icofcor}.  Thus, $J\subseteq D(I,M)$ and $\dim R/D(I,M)\le 1$.
\qed 

Finally, we prove part (d) of Proposition \ref{bigprop}:

\begin{prop} \label{bigprop4} Let $R$ be an
 unramified regular local ring, $I$ an ideal of $R$, and
$M$ a faithful 
finitely generated $R$-module.  Suppose that $\dim R/I\le 3$ and $M$ satisfies $S_{d-3}$ where
$d=\dim R$.
Then $\dim R/D(I,M)\le 1$.
\end{prop} 

{\it Proof:} By Proposition \ref{bigprop3}, we may assume that $\dim R/I=3$.
By Lemma \ref{fg}, $\hiM$ is finitely generated for $i<d-3$.

Suppose first that $R/I$ is
equidimensional.  

{\it Claim:} $\supp H^{d-2}_I(M)$ is finite.

{\it Proof:} Let $Q\in \supp H^{d-2}_I(M)$.
Clearly, $\hgt Q\ge d-2$.  In fact, 
$\hgt Q\ge d-1$ by the Hartshorne-Lichtenbaum Vanishing Theorem (HLVT), as $\dim R_Q/I_Q>
0$.  If $\hgt Q=d-1$ then $Q\in \supp H^{d-2}_I(R)$ (since, again by HLVT, $H^{d-1}_I(R)_Q=0$).
Thus, $Q\in \ass H^{d-2}_I(R)$, since $Q$ is minimal in $\supp H^{d-2}_I(R)$.
As $\ass H^{d-2}_I(R)$ is a finite set (\cite{HS}, \cite{L1}, \cite{L2}),
the claim now follows.  

By the claim, we have that $\dim R/D^{d-2}(I,M)\le 1$.  Now let 
$J= D^{d-2}(I,M)\cap D^{d-1}(I,M)\cap D^d(I,M)$.  Then $\dim R/J\le 1$ and for all
$x\in J$, $\hiR_x$ is finitely generated for all $i\neq d-3$.  
Thus, $\hiR_x$ is $I_x$-cofinite for all $i\ge 0$ and
$x\in J$ by Proposition \ref{Icof}.  Hence, $J\subseteq D(I,M)$. 

Now suppose that $R/I$ is not equidimensional.  Assume
$I=\sqrt{I}$ and let $I=K\cap L$ where $R/K$ is equidimensional,
$\dim R/K=3$, $\dim R/L\le 2$, and $\dim R/(K+L)\le 1$.   Now, let $J=D(K,M)\cap D(L,M)
\cap (K+L)$.
Then $\dim R/J\le 1$ by the preceding argument and Proposition \ref{bigprop3}.  
For all $x\in J$ we have that $K_x+L_x=R_x$ and
$H^i_I(M)_x\cong H^i_K(M)_x\oplus H^i_L(M)_x$ for all $i$.  Proceeding as in Case 2 of the proof of Proposition
\ref{bigprop2}, we obtain that $\hiM_x$ is $I_x$-cofinite for all $i$ and $x\in J$.  Thus,
$J\subseteq D(I,M)$ and $\dim R/D(I,M)\le 1$.
\qed

\section{Cofiniteness and Associated Primes}

In this section we apply the results of section 2 
to show that, under the conditions of
Theorem \ref{bigprop}, the set of associated primes of $H^i_I(M)$ is finite.  
Moreover, we prove that
the set $\ass \exti(N,\hjM)$ is finite for all $i$, $j$, 
and any finitely generated $R$-module $N$ with support in 
$\operatorname{V}(I)$. 
We first recall some facts from \cite{Mar}:

\begin{rmk} \label{complete} Let $(R,m)$ be a a local ring, $I$ an ideal of $R$, and $M$ an
$R$-module.  Let $\hat R$ be the $m$-adic completion of $R$.
\begin{enumerate}[(a)]
\item  If $\operatorname{Ass}_{\hat R}(M\otimes_R\hat R)$ is finite then $\ass M$ is finite.
\item  If $\supp M\subseteq \operatorname{V}(I)$ then $\ass M=\ass \operatorname{Hom}_R(R/I,M)$.
\end{enumerate}
\end{rmk}

We also need a slight improvement of a result due to Burch \cite{Bu}: 
 
\begin{lem}\label{Bulem}
Suppose $(R,m)$ is a complete local ring and $\{P_i\}_{i \in \mathbb{Z}^+}$ is a countable collection of prime ideals of $R$, none of which contain $I$.  Then there is an element $x \in I$ which is not contained in $P_i$ for all $i$.
\end{lem}
  
{\it Proof:}  Without loss of generality, we can assume there are no containment relations among the $P_i$.  We inductively construct a Cauchy sequence $\{x_1, x_2, x_3, \dots \} \subseteq I$ 
as follows:
Choose $x_1 \in I$ with $x_1\not\in P_1$.  
Suppose we have $x_1, x_2, \ldots, x_{r-1} \in I$ such that for all $i,s$ with 
$i \leq s \leq r-1$ we have $x_s \notin P_i$ and $x_s-x_i \in P_i\cap I^i$.

If $x_{r-1} \notin P_r$, set $x_r=x_{r-1}$.  If $x_{r-1} \in P_r$, 
choose $y_r$ such that $y_r \notin P_r$ but $y_r \in P_1 \cap P_2 \cap \cdots \cap P_{r-1} \cap I$.  Set 
$x_r=x_{r-1}+y_r^{r-1}$.  Note that $x_r \in I$ but $x_r \notin P_i$ for all $i \leq r$.  Clearly,
$$x_r-x_{r-1}=y_r^{r-1} \in P_1 \cap P_2 \cap \cdots \cap P_{r-1} \cap I^{r-1}.$$ 
Furthermore, for all $i\le r-1$
$$x_r-x_i=x_r-x_{r-1}+x_{r-1}-x_i \in  P_{i} \cap I^i.$$  
Assume that $x_r$ is chosen in this fashion for all positive integers $r$. Then $\{x_r\}_{r=1}^{\infty}$ is a Cauchy sequence in the complete local ring $R$ and hence has a limit $x$.  As
ideals are closed in the $m$-adic topology, $x\in I$.  Also, for any fixed $i\ge 1$, $\{x_s-x_i\}_{s=i}^{\infty}$ is a Cauchy sequence contained in $P_i$.  Thus $x-x_i\in P_i$ for all $i$ and hence
$x\not\in P_i$ for all $i$.
\qed

\begin{thm} \label{finass}  Let $R$ be a local ring, $I$ an ideal of $R$, and
$M$ a finitely generated $R$-module of dimension $n$.  Suppose one of the following conditions holds:
\begin{enumerate}[(a)]
\item $\dim M\le 3$;
\item $\dim R=4$ and the completion of $R$ is a UFD;
\item $R$ is the quotient of a Cohen-Macaulay ring, $\dim R/I\le 2$, and either $\dim M\le 4$ or 
$M$ satisfies Serre's condition $S_{n-3}$;
\item $R$ is a unramified regular local ring, $\dim R/I\le 3$, and $M$ satisfies $S_{d-3}$ where
$d=\dim R=\dim M$.
\end{enumerate}
Then for any finitely generated $R$-module $N$ such that $\supp N\subseteq \operatorname{V}(I)$
we have $\ass \exti (N,\hjM)$ is a finite set for all $i$ and $j$.  In particular, $\ass \hiM$
is finite for all $i$.
\end{thm}

{\it Proof:}  By Remark \ref{complete} we can assume that $R$ is complete.
Suppose that $\ass \exti (N,\hjM)$ is infinite for some $i$, $j$, $N$.
By Theorem \ref{bigprop}, $\dim R/D\le 1$ where
$D=D(I,M)$.  Hence,
there exists a countably infinite subset 
$\{P_{\ell}\}$
of $\ass \exti (N,\hjM)$ such that $P_{\ell}\not\supset D$ for all $\ell$.
By Lemma \ref{Bulem}, there exists $x\in D$ such that $x\not\in
P_{\ell}$ for all $\ell$.   As $x\in D$, 
$\ext^p_R (R/I,H^q_I(M))_x$ is a finitely generated $R_x$-module for
all $p,q$.  Hence, by \cite[Lemma 4.2]{HK} we have that $\exti (N,\hjM)_x$ is a finitely generated
$R_x$-module.  But as $x\not\in P_i$ for all $i$, 
$\operatorname{Ass}_{R_x} \exti (N,\hjM)_x$ is an infinite set, a contradiction.

The last assertion follows by Remark \ref{complete}.
\qed

In \cite{Kz}, M. Katzman proves that if $R=k[x,y,u,v,s,t]/(f)$, where $k$ is an arbitrary
field and $f=sx^2v^2-(t+s)xyuv+ty^2u^2$, and $I=(u,v)R$, then $H^2_I(R)$ has infinitely
many associated primes.  This shows that the conclusion of Theorem \ref{finass} does not hold in
general if $\dim R=5$ or $\dim R/I=4$, even when $M$ is Cohen-Macaulay.

Under the conditions of Theorem
\ref{bigprop}, we can also show that the non-$I$-cofinite locus of
$\hiM$ is finite for all $i$:

\begin{thm}\label{almostcof} Let $(R,m)$ be a local ring, $I$ an ideal of $R$, and
$M$ a finitely generated $R$-module of dimension $n$.  Suppose one of the following conditions holds:
\begin{enumerate}[(a)]
\item $\dim M\le 3$;
\item $\dim R/I\le 2$;
\item $\dim R=4$ and $R$ is a UFD;
\item $R$ is an unramified regular local ring, $\dim R/I\le 3$, and $M$ satisfies $S_{d-3}$
where $d=\dim R=\dim M$.
\end{enumerate}
Then $\hiM_p$ is $I_p$-cofinite for all but finitely many primes $p$ of $R$.
Furthermore, the Bass numbers $\mu_j(p,\hiM)$ are finite for $i$, $j$ and all but finitely many
primes $p$.
\end{thm}

{\it Proof:} The second assertion follows readily from the first and \cite[Lemma 4.2]{HK}.
Case (b) follows from \cite[Theorem 1]{DM}, since $\dim R_Q/I_Q\le 1$ for all $Q\neq m$.
For cases (a), (c), and (d), we have that $\dim R/D\le 1$ where $D=D(I,M)$ by Theorem 
\ref{bigprop}.  Therefore, $\operatorname{V}(D)$ is a finite set and for all $p\not\in 
\operatorname{V}(D)$, $\hiM_p$ is $I_p$-cofinite for all $i$. 
\qed

We now wish to give some examples to show
that the conclusion of Theorem \ref{almostcof} (and therefore of Theorem \ref{bigprop}) does
not hold in general. Before doing so, we first prove
the following result, which is the core argument in the proof of \cite[Theorem 2.3]{HK}:

\begin{lem} \label{hom} Let $(R,m)$ be a complete local ring, $I$ an ideal of $R$, and
$N$ an $R$-module such that $\supp N\subseteq \{m\}$.  Suppose $\Hom_R(R/I,N)$ is
finitely generated.  Then $N$ is Artinian, 
$I+\ann_RN$ is $m$-primary, and $\ext^i_R(R/I,N)$ has finite
length for all $i$.
\end{lem}

{\it Proof:} Since $\Hom_R(R/m,N)$ is isomorphic to a submodule of $\Hom_R(R/I,N)$,
we have that $\Hom_R(R/m,N)$ is finitely generated.  Consequently, as $\supp N\subseteq
\{m\}$, $N$ is Artinian and the Matlis dual $N^{\vee}$ of $N$ is a finitely generated $R$-module.
Now, by \cite[Theorem 11.57]{Ro}, 
$\ext^i_R(R/I,N)^{\vee}\cong
\operatorname{Tor}^R_i(R/I,N^{\vee})$. 
As $\Hom_R(R/I,N)$
is finitely generated and Artinian, it has finite length.  Consequently, 
$R/I\otimes_R N^{\vee}$
has finite length and $I+\ann_R(N^{\vee})$ is $m$-primary.  This implies that
$\operatorname{Tor}^R_i(R/I,N^{\vee})$ has finite length for all $i$.  The last assertion now
follows by Matlis duality. 
\qed

\begin{prop} \label{notIcof} Let $(R,m)$ be a 
$d$-dimensional analytically normal local Cohen-Macaulay
domain and
$I$ an ideal of $R$.  Suppose
\begin{enumerate}[(a)]
\item $\dim R/I\ge 2$, and
\item $\dim R/Q=1$ for some $Q\in \operatorname{Min}_RR/I$.
\end{enumerate}
Then $\Hom_R(R/I,H^{d-1}_I(R))$ 
is not finitely generated.
\end{prop}

{\it Proof:} By passing to the completion, we can assume $R$ is a complete CM normal domain.
Suppose $\Hom_R(R/I,H^{d-1}_I(R))$ is finitely generated.  Assume that
$I=\sqrt{I}$ and let $I=J\cap K$ where
$\dim R/J=1$, $\dim R/p\ge 2$ for all $p\in \operatorname{Min}_RR/K$, and $\sqrt{J+K}=m$.
By HLVT, we have
$H^d_J(R)=H^d_K(R)=0$.  Also, $\supp H^{d-1}_K(R)\subseteq \{m\}$, 
since for all dimension one primes $P$ containing
$K$ one has $\dim (R/K)_P>0$ and $R_P$ is analytically irreducible.
By the Mayer-Vietoris
sequence we have
$$0\to H^{d-1}_J(R)\oplus H^{d-1}_K(R)\to H^{d-1}_I(R)\to H^d_m(R)\to 0.$$ 
Thus, we have the long exact sequence
$$0\to \Hom_R(R/J,H^{d-1}_J(R))\oplus \Hom_R(R/J,H^{d-1}_K(R))\to \Hom_R(R/J,H^{d-1}_I(R))\to$$
$$\Hom_R(R/J,H^d_m(R))
\to \ext^1_R(R/J,H^{d-1}_J(R))\oplus \ext^1_R(R/J,H^{d-1}_K(R))\to \cdots$$
Now $\Hom_R(R/J,H^{d-1}_I(R))$ is finitely generated as it is isomorphic to a submodule of
$\Hom_R(R/I,H^{d-1}_I(R))$.  Also, since
$\dim R/J=1$, $\ext^i_R(R/J,H^{d-1}_J(R))$ is finitely generated for all $i$ (\cite[Theorem 1]{DM}).
Further, as $\Hom_R(R/J,H^{d-1}_K(R))$ is finitely generated (by the long exact sequence
above) and $\supp H^{d-1}_K(R)\subseteq \{m\}$, we can apply Lemma \ref{hom} to get that 
$\ext^i_R(R/J,H^{d-1}_K(R))$ has finite length for all $i$.
Hence, $\Hom_R(R/J,H^d_m(R))$
is finitely generated.  Applying Lemma \ref{hom} again, we have $J+\ann_RH^d_m(R)$ is
$m$-primary, which is a contradiction since $\ann_RH^d_m(R)=0$. 
\qed

\begin{ex} {\rm Let $R=k[x,y,z]_{(x,y,z)}$ where $k$ is a field.  
Let $I=((x)\cap (y,z))R$.  Then $\Hom_R(R/I,H^2_I(R))$
is not finitely generated.   Consequently, $H^i_I(R)$ is not $I$-cofinite for $i=1,2$
by Proposition \ref{Icof}.}
\end{ex}

In a sense, this example represents the ``minimal'' example of a
local cohomology module which is not cofinite.  If $\dim R\le 2$ then $H^i_I(M)$
is $I$-cofinite for all $I$, $M$ and $i\ge 0$ by Corollary \ref{Icofcor}(a).
If $R$ is a $3$-dimensional regular local ring and
$I$ is an ideal such that $H^i_I(R)$ is not $I$-cofinite for
some $i$, then one can easily show that $I=(f)\cap J$ where $\dim R/J=1$ and $(J,f)$
is $m$-primary.

\begin{ex} {\rm Let $R=k[x,y,z,u,v]_{(x,y,z,u,v)}$, $I=((x)\cap(y,z))R$, and $P=(x,y,z)R$.  Then
$\Hom_R(R/I,H^2_I(R))_P$ is not a finitely generated $R_P$-module by Proposition \ref{notIcof}.
Hence  
$\Hom_R(R/I,H^2_I(R))_Q$
is not a finitely generated $R_Q$-module for any prime $Q$ of $R$ containing $P$.}
\end{ex}  

The above example shows there exists
a five-dimensional regular local ring $R$ and an ideal $I$ with
$\dim R/I=4$ such that $H^i_I(R)_Q$ is not $I_Q$-cofinite for $i=1,2$ and
infinitely many primes $Q$.  Hence, $\dim R/D(I,R)\ge 2$. 
However, the Bass numbers and the sets of associated primes
of $H^i_I(R)$ are finite for all $i$ for this example by \cite{HS} and \cite{L1}.

Our final result extends Proposition \ref{notIcof} and generalizes \cite[Theorem 3.6(ii)]{HK}.  

\begin{thm} \label{notIcof2} Let $(R,m)$ be an analytically normal Cohen-Macaulay
local domain and $I$
an ideal such that $\dim R/I\ge 2$.  If $\operatorname{Spec} R/I-\{m/I\}$ is
disconnected then $\operatorname{Hom}_R(R/I,H^{d-1}_I(R))$ is not finitely generated.
\end{thm}

{\it Proof:} We may assume $R$ is a complete CM normal domain. By
virtue of Proposition \ref{notIcof}, we
can assume $\dim R/Q\ge 2$ for all $Q\in \operatorname{Min}_RR/I$.
Hence, $\supp H^{d-1}_I(R)\subseteq \{m\}$, since for all primes
$P$ of height $d-1$ containing $I$ we have $\dim (R/I)_P>0$ and $R_P$ is analytically
irreducible.  Now, as the punctured spectrum of $R/I$ is disconnected, we can use
the Mayer-Vietoris sequence to obtain a surjective map $H^{d-1}_I(R)\to H^d_m(R)$.
Since $\ann_RH^d_m(R)=0$ we obtain that $\ann_RH^{d-1}_I(R)=0$.
Hence $\operatorname{Hom}_{R}(R/I,H^{d-1}_I(R))$ is not finitely generated 
by Lemma \ref{hom}.

\begin{ex} {\rm Let $k$ be a field, $R=(k[x,y,u,v]/(xu-yv))_{(x,y,u,v)}$, 
and $I=(x,y)R\cap (u,v)R$.  Then
$\Hom_R(R/I,H^2_I(R))$ is not finitely generated.}
\end{ex}


\begin{thebibliography}{99}


\bibitem[BS]{BS}
Brodmann,~M. and Sharp,~R., {\em Local Cohomology: an algebraic introduction with geometric applications}, Cambridge Studies in Advanced Mathematics no. {\bf 60}, Cambridge, Cambridge University Press, 1998.

\bibitem[BH]{BH}
Bruns,~W. and Herzog,~J., {\em Cohen-Macaulay Rings}, Cambridge Studies in Advanced Mathematics no. {\bf 39}, Cambridge, Cambridge University Press, 1993.

\bibitem[Bu]{Bu}
Burch,~L., {\em Codimension and analytic spread}, Proc. Camb. Phil. Soc. {\bf 72}, 369-373 (1972).

\bibitem[DM]{DM}
Delfino,~D. and Marley,~T., {\em Cofinite modules and local cohomology}, J. Pure and App. Alg. {\bf 121}, 45-52 (1997).

\bibitem[Fa]{Fa}
Faltings,~G., {\em \"Uber die Annulatoren lokaler Kohomologiegruppen}, Archiv der Math. {\bf 30}, 
473-476 (1978).

\bibitem[EGA]{EGA}
Grothendieck,~A., {\em \'El\'ements de g\'eom\'etrie alg\'ebrique}, Inst. Hautes \'Etudes Sci.
Publ. Math. {\bf 24} (1965).

\bibitem[Ha]{Ha}
Hartshorne,~R., {\em Affine Duality and Cofiniteness}, Inv. Math. {\bf 9}, 145-164 (1970).


\bibitem[HK]{HK}
Huneke,~C. and Koh, J., {\em Cofiniteness and vanishing of local cohomology modules}, Math. Proc Camb. Phil. Soc. {\bf 110}, 421-429 (1991). 

\bibitem[HS]{HS}
Huneke,~C. and Sharp,~R., {\em Bass Numbers of local cohomology modules}, Trans. A.M.S. {\bf 339}, 765-779 (1993).

\bibitem[Kz]{Kz}
Katzman,~M., {\em An example of an infinite set of associated primes of a local cohomology module}, J. Algebra, to appear.

\bibitem[Ka]{Ka}
Kawasaki,~K.-I., {\em Cofiniteness of local cohomology modules for principal ideals}, Bull. London
Math. Soc. {\bf 30}, 241-246 (1998).



\bibitem[L1]{L1}
Lyubeznik,~G., {\em Finiteness Properties of local cohomology modules (An application of $D$-modules to commutative algebra)}, Inv. Math. {\bf 113}, 41-55 (1993).

\bibitem[L2]{L2}
Lyubeznik,~G., {\em Finiteness properties of local cohomology modules for regular local rings of mixed characteristic: The unramified case}, Comm. Alg. {\bf 28} no. 12, 5867-5882 (2000).

\bibitem[Mar]{Mar}
Marley,~T., {\em The associated primes of local cohomology modules of small dimension}, Manuscripta Math. {\bf 104}, 519-525 (2001).

\bibitem[Mat]{Mat}
Matsumura,~H., {\em Commutative Ring Theory}, Cambridge Studies in Advanced Mathematics no. {\bf 8}, Cambridge, Cambridge University Press, 1986.

\bibitem[Ro]{Ro}
Rotman,~J., {\em An Introduction to Homological Algebra}, Orlando, FL, Academic Press, 1979.

\bibitem[Si]{Si}
Singh,~A., {\em $p$-torsion elements in local cohomology modules}, Math. Res. Letters {\bf 7}, 165-176 (2000).


\bibitem[Yo]{Yo}
Yoshida,~K.-I. {\em Cofinitness of local cohomology modules for ideals of dimension one}, Nagoya 
Math. J. {\bf 147}, 179-191 (1997).


\end{thebibliography}
\end{document}